\documentclass[12pt, reqno]{amsart}
\usepackage{amsmath, amsthm, amscd, amsfonts, amssymb}
\usepackage[bookmarksnumbered]{hyperref}

\newtheorem{theorem}{Theorem}[section]

\newtheorem{corollary}[theorem]{Corollary}
\theoremstyle{definition}

\theoremstyle{remark}

\numberwithin{equation}{section}

\begin{document}

\title[A Characterization of $(\sigma,\tau)-$ derivations  on von Neumann algebras ]
{A Characterization of $(\sigma,\tau)-$ derivations  on von Neumann
algebras }
\author[M. Eshaghi Gordji]{M. Eshaghi Gordji }
\address{Department of Mathematics, Semnan University, P. O. Box 35195-363, Semnan, Iran}
\email{madjid.eshaghi@gmail.com}



\begin{abstract} Let $\mathcal A$ be a von Neumann algebra and $\mathcal M$ be a Banach  $\mathcal A-$module.
It is shown that for every  homomorphisms $\sigma, \tau$ on
$\mathcal A$, every bounded linear map $f:\mathcal A\to \mathcal M$
 with property that $f(p^2)=\sigma(p)f(p)+f(p)\tau(p)$
for every projection $p$ in $\mathcal A$  is a
$(\sigma,\tau)-$derivation. Also, it is shown that a bounded linear
map $f:\mathcal A \to \mathcal M $ which satisfies $f(ab)=
\sigma(a)f(b)+f(a)\tau(b)$ for all $a,b\in \mathcal A$ with $ab=S$,
is a $(\sigma,\tau)-$ derivation if $\tau(S)$ is left invertible for
fixed $S$.

\end{abstract}
\maketitle


\section{Introduction}
Let $\mathcal A$ be a Banach algebra. An $\mathcal A-$module
$\mathcal M$ is a Banach $\mathcal A-$module if $\mathcal M$ is a
Banach space and the $\mathcal A-$module maps $(a; x) \to  ax;
\mathcal A \times \mathcal M \to  \mathcal M;$ and $(x; a) \to xa;
\mathcal M\times \mathcal A \to \mathcal M;$ satisfy
$max\{\|ax\|,\|xa\|\}\leq \|a\|\|x\|$ for all $a\in \mathcal A$ and
$x\in \mathcal M$.  Suppose that $\mathcal A$ is unital. We denote
the identity of $\mathcal A$ by 1. A Banach $\mathcal A-$module
$\mathcal M$ is called unital provided that $1 x = x = x1$ for each
$x\in \mathcal M$.

Recently, a number of authors \cite{ B-V, M-M1, M-M2} have studied
various generalized notions of derivations in the context of Banach
algebras. There are some applications in the other fields of
research \cite{H-L-S}. Such mappings have been extensively studied
in pure algebra; cf. \cite{A-R, B, H}. A generalized concept of
derivation is as follows.

Let $\mathcal A$ be a Banach algebra and $\mathcal M$ be  a Banach
$\mathcal A-$module. Let $\sigma,\tau \in \mathcal {BL}(\mathcal A)$
be bounded
 linear maps   on ${\mathcal  A}$.  A linear mapping $d:\mathcal A
\to\mathcal M$ is called a

$\bullet$ $(\sigma,\tau)-$derivation  if
$$d(ab)=\sigma(a)d(b)+d(a)\tau(b)\quad (a, b \in \mathcal A).\eqno (1)$$

$\bullet$ $(\sigma,\tau)-$Jordan derivation if
$$d(a^2)=\sigma(a)d(a)+d(a)\tau(a)\quad (a \in \mathcal A).\eqno(2)$$

For instance  every ordinary derivation (Jordan derivation) of an
algebra $\mathcal A$ into an $\mathcal A$--module $\mathcal M$ is an
$(id_ {\mathcal A},id_ {\mathcal A})-$derivation ($(id_ {\mathcal
A},id_ {\mathcal A})-$Jordan derivation), where $id_{\mathcal A}$ is
the identity mapping on the algebra $\mathcal A$. As another
example, every homomorphism (Jordan homomorphism) $h:\mathcal A \to
\mathcal A$ is a $(\frac{h}{2},\frac{h}{2})-$derivation
($(\frac{h}{2},\frac{h}{2})-$Jordan derivation).

Clearly, every  $(\sigma,\tau)-$derivation is a
$(\sigma,\tau)-$Jordan derivation. Using the fact that $ab+ba =
(a+b)^2-a^2-b^2$, it is easy to prove that the
$(\sigma,\tau)-$Jordan derivation identity is equivalent to
$$d(ab+ba)=\sigma(a)d(b)+d(a)\tau(b)+\sigma(b)d(a)+d(b)\tau(a) \quad (a, b \in \mathcal A). \eqno(3)$$

We refer to \cite{P} for the general theory of these notions.

\section{Main result}

 In 1996,
Johnson  \cite{Jo} proved the following theorem (see also Theorem
2.4 of \cite{H-L}).

\begin{theorem}
Suppose $\mathcal A$ is a $C^*-$ algebra and $\mathcal M$ is a
Banach $\mathcal A-$module.  Then each Jordan derivation
 $d:\mathcal A\to \mathcal M$ is a derivation.
\end{theorem}
As an application of this theorem, we give the following result for
characterization  of $(\sigma,\tau)-$derivations on von Neumann
algebras.
\begin{theorem} Let  $\mathcal A$ be a von Neumann algebra and let $\sigma, \tau$ be bounded
homomorphisms on $\mathcal A$. Let  $\mathcal M$ be a Banach
$\mathcal A-$module and  $d:\mathcal A\to \mathcal M$ be a bounded
linear map with property that $d(p^2)=\sigma(p)d(p)+d(p)\tau(p)$ for
every projection $p$ in $\mathcal A$. Then $d$   is a
$(\sigma,\tau)-$derivation.
\end{theorem}
\begin{proof}
We prove the theorem in two steps as follows. \\
STEP I.  Recall that $\sigma, \tau$ are bounded homomorphisms  on
$\mathcal A$ and $d:\mathcal A\to \mathcal M$ is a bounded linear
map with property that $d(p^2)=\sigma(p)d(p)+d(p)\tau(p)$ for every
projection $p$ in $\mathcal A$. We show that $d$ is a Jordan
homomorphism. Let $p,q\in \mathcal A$ be orthogonal projections in
$\mathcal A$. Then $p+q$ is a projection wherefore by assumption,
\begin{align*}
\sigma(p)d(p)&+d(p)\tau(p)+
\sigma(q)d(q)+d(q)\tau(q)=d(p)+d(q)=d(p+q)\\
&=\sigma(p+q)d(p+q)+d(p+q)\tau(p+q)=
\sigma(p)d(p)+d(p)\tau(p)\\
&+\sigma(q)d(q) +d(q)\tau(q)+\sigma(p)d(q)+d(p)\tau(q)
+\sigma(q)d(p)+d(q)\tau(p).
\end{align*}
This means that  $$\sigma(p)d(q)+d(p)\tau(q)
+\sigma(q)d(p)+d(q)\tau(p)=0. \eqno (4)$$ Let
$a=\sum_{j=1}^{n}\lambda_jp_j$ be a combination of mutually
orthogonal projections $p_1,p_2,...,p_n \in \mathcal A.$ Then we
have $$\sigma(p_i)d(p_j)+d(p_i)\tau(p_j)
+\sigma(p_j)d(p_i)+d(p_j)\tau(p_i)=0 \eqno (5)$$ for all $i,j\in
\{1,2,...,n\}$ with $i\neq j.$ So
$$d(a^2)=d(\sum_{j=1}^{n}\lambda_j^2p_j)=\sum_{j=1}^{n}\lambda_j^2d(p_j).\eqno (6)$$
On the other hand by (5),  we obtain that
\begin{align*}
\sigma(a)d(a)&+d(a)\tau(a)=
\sigma(\sum_{j=1}^{n}\lambda_jp_j)\sum_{j=1}^{n}\lambda_jd(p_j)
+\sum_{j=1}^{n}\lambda_jd(p_j)\tau(\sum_{j=1}^{n}\lambda_jp_j)\\
&=\lambda_1d(p_1)\tau(\sum_{j=1}^{n}\lambda_jp_j)+\lambda_2d(p_2)\tau(\sum_{j=1}^{n}\lambda_jp_j)
+...+\lambda_nd(p_n)\tau(\sum_{j=1}^{n}\lambda_jp_j)\\
&+\sigma(\lambda_1p_1)\sum_{j=1}^{n}\lambda_jd(p_j)+\sigma(\lambda_2p_2)\sum_{j=1}^{n}\lambda_jd(p_j)+
...+\sigma(\lambda_np_n)\sum_{j=1}^{n}\lambda_jd(p_j)\\
&=\lambda_1^2(\sigma(p_1)d(p_1)+d(p_1)\tau(p_1))+\lambda_2^2(\sigma(p_2)d(p_2)+d(p_2)\tau(p_2))+...\\
&+
\lambda_n^2(\sigma(p_n)d(p_n)+d(p_n)\tau(p_n))=\sum_{j=1}^n\lambda_j^2d(p_j).\hspace{5.7cm}(7)
\end{align*}
Combining (6) by (7) to get $d(a^2)=\sigma(a)d(a)+d(a)\tau(a)$. By
the spectral theorem (see Theorem 5.2.2 of  \cite{Ka-Ri}), every
self adjoint element $a\in \mathcal A_{sa}$ is the norm--limit of
finite combinations of mutually orthogonal projections. Since
$d,\sigma,\tau$ are bounded, then
$$d(a^2)=\sigma(a)d(a)+d(a)\tau(a) \eqno(8)$$ for all $a\in \mathcal A_{sa}$.
Replacing $a$ by $a+b$ in (8), we obtain
$$d(ab+ba)=\sigma(a)d(b)+d(a)\tau(b)+\sigma(b)d(a)+d(b)\tau(b)\eqno (9)$$
for all $a,b\in \mathcal A_{sa}$. Let $a\in \mathcal A$. Then there
are $a_1,a_2\in \mathcal A_{sa}$ such that $a=a_1+ia_2$. Hence,
\begin{align*}d(a^2)&=d(a_1^2+a_2^2+i(a_1a_2+a_2a_1))=d(a_1^2)+d(a_2^2)+id(a_1a_2+a_2a_1)\\
&=\sigma(a_1)d(a_1)+d(a_1)\tau(a_1)+\sigma(a_2)d(a_2)+d(a_2)\tau(a_2)\\
&+i(\sigma(a_1)d(a_2)+d(a_1)\tau(a_2)+\sigma(a_2)d(a_1)+d(a_2)\tau(a_1))\\
&=(\sigma(a_1)+i\sigma(a_2))(d(a_1)+id(a_2))+(d(a_1)+id(a_2))(\tau(a_1)+i\tau(a_2))\\
&=\sigma(a)d(a)+d(a)\tau(a).
\end{align*}
STEP II. We show that every $(\sigma,\tau)-$ Jordan derivation from
$\mathcal A$ into $\mathcal M$ is a $(\sigma,\tau)-$derivation. Let
$d:\mathcal A \to \mathcal M$ be a $(\sigma,\tau)-$ Jordan
derivation. It is easy to see that $\mathcal M$ is a Banach
$\mathcal A-$module by the following module actions:
\begin{eqnarray*}
a\cdot m=\sigma(a)m,\qquad m\cdot a=m\tau(a)\qquad (a\in\mathcal A,
m\in\mathcal M)
\end{eqnarray*}
we denote $\mathcal M_{(\sigma,\tau)}$ the above $\mathcal
A-$module. One can show that $d$ is a Jordan derivation from
$\mathcal A$ into $\mathcal M_{(\sigma,\tau)}$. It follows form
Theorem 2.1 that $d$ is a derivation from $\mathcal A$ into
$\mathcal M_{(\sigma,\tau)}$. Hence,  $d$ is a $(\sigma,\tau)-$
derivation from $\mathcal A$ into $\mathcal M$.
\end{proof}

Suppose that $\mathcal A$ is a Banach algebra and $\mathcal M$ is an
$\mathcal A-$module. Let $S$ be in A. We say that $S$ is right
separating point of $\mathcal M$ if the condition  $mS = 0$ for
$m\in \mathcal M$ implies $m = 0$.

\begin{theorem}   Let $\mathcal A$ be a unital Banach algebra and $\mathcal M$ be a Banach  $\mathcal A-$module.
Let $S$  be in $\mathcal A$ and  $ \tau\in Hom(\mathcal A)$ be
bounded homomorphism on $\mathcal A$ which satisfies $\tau(S)$ is a
right separating point of $\mathcal M$ and let $\sigma$ be a bounded
Jordan homomorphism on $\mathcal A$. Let   $f:\mathcal A \to
\mathcal M $ be a bounded linear map. Then the following assertions
are equivalent

a) $f(ab)= \sigma(a)f(b)+f(a)\tau(b) $ for all $a,b\in \mathcal A$
with $ab=S$.

b) $f$ is a $(\sigma,\tau)-$Jordan derivation which satisfies
 $f(Sa)=
\sigma(S)f(a)+f(S)\tau(a)$  and $f(aS)= \sigma(a)f(S)+f(a)\tau(S)$
for all $a\in \mathcal A.$
\end{theorem}
\begin{proof}
First suppose that (a) holds. Then we have
$$f(S)=f(1S)=\sigma(1)f(S)+f(1)\tau(S)=f(S)+f(1)\tau(S)$$
hence, by assumption, we get that $f(1)=0.$ Let $a\in \mathcal A.$
For scalars $\lambda$ with $|\lambda| < \frac {1}{\|a\|},$
$1-\lambda a$ is invertible in  $\mathcal A$. Indeed, $(1-\lambda
a)^{-1}=\sum_{n=0}^{\infty} \lambda^na^n$. Then
\begin{align*} f(S)&=f[(1-\lambda a)(1-\lambda a)^{-1}S]=\sigma ((1-\lambda
a))f((1-\lambda a)^{-1}S)\\
&+f((1-\lambda a))\tau((1-\lambda a)^{-1}S)=\sigma ((1-\lambda
a))f(\sum_{n=0}^{\infty} \lambda^na^nS)\\
&-\lambda f(a)\tau(\sum_{n=0}^{\infty}
\lambda^na^nS)=f(S)+\sum_{n=1}^{\infty}
\lambda^n[f(a^nS)\\
&-f(a)\tau(a^{n-1}S)-\sigma(a)f(a^{n-1}S)].
\end{align*}
So $$\sum_{n=1}^{\infty}
\lambda^n[f(a^nS)-f(a)\tau(a^{n-1}S)-\sigma(a)f(a^{n-1}S)]=0$$ for
all $\lambda$ with $|\lambda| < \frac {1}{\|a\|}.$ Consequently
$$f(a^nS)-f(a)\tau(a^{n-1}S)-\sigma(a)f(a^{n-1}S)=0 \eqno (10)$$
for all $n\in \Bbb N$. Put $n=1$ in (10) to get
$$f(aS)=\sigma(a)f(S)+f(a)\tau(S). \eqno (11)$$
Similarly, using equation  $f(S)=f[S(1-\lambda a)^{-1}(1-\lambda
a)]$ we get
$$f(Sa)=\sigma(S)f(a)+f(S)\tau(a)$$
for all $a\in \mathcal A.$

Now, put $n=2$ in (10) to get
$$f(a^2S)=\sigma(a)f(aS)+f(a)\tau(aS). \eqno (12)$$
Combining (11), (12) to obtain
$$f(a^2S)=\sigma(a)(\sigma(a)f(S)+f(a)\tau(S))+f(a)\tau(aS). \eqno (13)$$
Replacing $a$ by $a^2$ in (11), we get
$$f(a^2S)=\sigma(a^2)f(S)+f(a^2)\tau(S). \eqno (14)$$
It follows from (13), (14) that
$$(f(a^2)-\sigma(a)f(a)-f(a)\tau(a))\tau(S)=0.\eqno(15)$$
On the other hand $\tau(S)$ is right separating point of $\mathcal
M.$ Then by (15) $f$ is a $(\sigma,\tau)-$Jordan derivation.

Now suppose that the condition (b) holds. Let $a,b\in \mathcal A$
which satisfy $ab=S.$ The next relation  follows from a
straightforward computation using the $(\sigma,\tau)-$Jordan
derivation identities (2) and (3).
\begin{align*}
f(Sa)&=f(aba)=\frac{1}{2}[f(a(ab+ba)+(ab+ba)a)-f(a^2b+ba^2)]\\
&=\frac{1}{2}[f(a)\tau(ab+ba)+\sigma(a)f(ab+ba)+f(ab+ba)\tau(a)+\sigma(ab+ba)f(a)\\
&-f(a^2)\tau(b)-\sigma(a^2)f(b)-f(b)\tau(a^2)-\sigma(b)f(a^2)]\\
&=f(a)\tau(ba)+\sigma(a)f(b)\tau(a)+\sigma(ab)f(a)\\
&=f(a)\tau(ba)+\sigma(a)f(b)\tau(a)+f(Sa)-f(S)\tau(a).
\end{align*}
So $$[f(S)-f(a)\tau(b)-\sigma (a)f(b)]\tau(a)=0.$$  Hence,
$$[f(S)-f(a)\tau(b)-\sigma
(a)f(b)]\tau(a)\tau(b)=[f(S)-f(a)\tau(b)-\sigma (a)f(b)]\tau(S)=0.$$
Since $\tau(S)$ is a right separating point of $\mathcal M,$ then
$$f(S)=f(a)\tau(b)+\sigma (a)f(b).$$
\end{proof}

By  Theorems 2.2 and 2.3, we have the following corollaries.
\begin{corollary} Let  $\mathcal A$ be a von Neumann algebra and let $\sigma, \tau$ be bounded
homomorphisms on $\mathcal A$. Let  $\mathcal M$ be a Banach
$\mathcal A-$module and  $d:\mathcal A\to \mathcal M$ be a bounded
linear map. Then the following assertions are equivalent

a) $d(p^2)=\sigma(p)d(p)+d(p)\tau(p)$ for every projection $p$ in
$\mathcal A$.

b) $\sigma(a)d({a}^{-1})+d(a)\tau({a}^{-1})=0$ for all invertible
$a\in \mathcal A.$

c)  $d$   is a $(\sigma,\tau)-$derivation.
\end{corollary}

\begin{corollary} Let  $\mathcal A$ be a von Neumann algebra and let   $\mathcal M$ be a Banach
$\mathcal A-$module and  $d:\mathcal A\to \mathcal M$ be a bounded
linear map. Then the following assertions are equivalent

a) $d(p^2)=pd(p)+d(p)p$ for every projection $p$ in $\mathcal A$.

b) $ad({a}^{-1})+d(a){a}^{-1}=0$ for all invertible $a\in \mathcal
A.$

c)  $d$   is a derivation.
\end{corollary}

\end{document}